\documentclass[12pt]{article}%
\usepackage{graphicx}
\usepackage[intlimits]{amsmath}
\usepackage{latexsym}
\usepackage{amsfonts}
\usepackage{amssymb}%
 \makeatletter
\newfont{\footsc}{cmcsc10 at 8truept}
\newfont{\footbf}{cmbx10 at 8truept}
\newfont{\footrm}{cmr10 at 10truept}
\newtheorem{theorem}{Theorem}

\newtheorem{claim}[theorem]{Claim}

\newenvironment{proof}[1][Proof]{\noindent{\textbf {#1}  }}  {\hfill$\Box$\bigskip}

\begin{document}

\title{Complete $r$-partite subgraphs of dense $r$-graphs}
\author{Vladimir Nikiforov\\{\small Department of Mathematical Sciences, University of Memphis, Memphis TN
38152}}
\maketitle

\begin{abstract}
Let $r\geq3$ and $\left(  \ln n\right)  ^{-1/\left(  r-1\right)  }\leq
\alpha\leq r^{-3}.$ We show that:

Every $r$-uniform graph on $n$ vertices with at least $\alpha n^{r}/r!$ edges
contains a complete $r$-partite graph with $r-1$ parts of size $\left\lfloor
\alpha\left(  \ln n\right)  ^{1/\left(  r-1\right)  }\right\rfloor $ and one
part of size $\left\lceil n^{1-\alpha^{r-2}}\right\rceil .$

This result follows from a more general digraph version:

Let $U_{1},\ldots,U_{r}$ be sets of size $n,$ and $M\subset U_{1}\times
\cdots\times U_{r}$ satisfy $\left\vert M\right\vert \geq\alpha n^{r}.$ If the
integers $s_{1},\ldots,s_{r-1}$ satisfy $1\leq s_{1}\cdots s_{r-1}%
\leq\left\lfloor \alpha^{r-1}\ln n\right\rfloor ,$ then there exists
$V_{1}\times\cdots\times V_{r}\subset M,$ such that $V_{i}\subset U_{i}$ and
$\left\vert V_{i}\right\vert =s_{i}$ for $1\leq i<r,$ and $\left\vert
V_{r}\right\vert >n^{1-\alpha^{r-2}}.$\medskip


\end{abstract}

In this note \emph{graph} means $r$\emph{-uniform graph} for some fixed
$r\geq3$.

Given $c>0,$ how large complete $r$-partite graphs\ must contain a graph $G$
with $n$ vertices and $cn^{r}$ edges? This question was answered for $r=2$ in
\cite{ErSt46}, and for $r>2$ in \cite{Erd64}: $G$ contains a complete
$r$-partite graph with each part of size $a\left(  \log n\right)  ^{1/\left(
r-1\right)  }$ for some $a=a\left(  c\right)  >0,$ independent of $n$.

Here we refine this statement for $r\geq3$ and extend it to digraphs. Letting
$K_{r}\left(  s_{1},\ldots,s_{r}\right)  $ be the complete $r$-partite graph
with parts of size $s_{1},\ldots,s_{r},$ our most concise result reads as:

\begin{theorem}
\label{th1}Let $r\geq3$ and $\left(  \ln n\right)  ^{-1/\left(  r-1\right)
}\leq\alpha\leq r^{-3}.$ Every graph with $n$ vertices and at least $\alpha
n^{r}/r!$ edges contains a $K_{r}\left(  s,\ldots,s,t\right)  $ with
$s=\left\lfloor \alpha\left(  \ln n\right)  ^{1/\left(  r-1\right)
}\right\rfloor $ and $t=\left\lceil n^{1-\alpha^{r-2}}\right\rceil .$
\end{theorem}

Theorem \ref{th1} follows immediately from a subtler one:

\begin{theorem}
\label{th2}Let $r\geq3$ and $\left(  \ln n\right)  ^{-1/\left(  r-1\right)
}\leq\alpha\leq r^{-3}$. Let $G$ be a graph with $n$ vertices and at least
$\alpha n^{r}/r!$ edges. If the integers $s_{1},\ldots,s_{r-1}$ satisfy $1\leq
s_{1}\cdots s_{r-1}\leq\alpha^{r-1}\ln n,$ then $G$ contains a $K_{r}\left(
s_{1},\ldots,s_{r-1},t\right)  $ with $t>n^{1-\alpha^{r-2}}.$
\end{theorem}

It seems that a digraph setup is more natural for such results, e.g., Theorem
\ref{th2} follows from

\begin{theorem}
\label{th3}Let $r\geq3$ and $\left(  \ln n\right)  ^{-1/\left(  r-1\right)
}\leq\alpha\leq r^{-3}.$ Let $U_{1},\ldots,U_{r}$ be sets of size $n$ and
$M\subset U_{1}\times\cdots\times U_{r}$ satisfy $\left\vert M\right\vert
\geq\alpha n^{r}.$ If the integers $s_{1},\ldots,s_{r-1}$ satisfy $1\leq
s_{1}\cdots s_{r-1}\leq\alpha^{r-1}\ln n,$ then there exists $V_{1}%
\times\cdots\times V_{r}\subset M$ such that $V_{i}\subset U_{i}$ and
$\left\vert V_{i}\right\vert =s_{i}$ for $1\leq i<r,$ and $\left\vert
V_{r}\right\vert >n^{1-\alpha^{r-2}}.$
\end{theorem}

We prove Theorem \ref{th3} by counting. For a better view on the matter we
give a separate theorem, hoping that it may have other applications as well.

Let $U_{1},\ldots,U_{r}$ be nonempty sets and $M\subset U_{1}\times
\cdots\times U_{r},$ let the positive integers $s_{1},\ldots,s_{r}$ satisfy
$\left\vert U_{i}\right\vert \geq s_{i}$ for $1\leq i\leq r$. Write
$B_{M}\left(  s_{1},\ldots,s_{r}\right)  $ for the set of products
$V_{1}\times\cdots\times V_{r}\subset M$ such that $V_{i}\subset U_{i}$ and
$\left\vert V_{i}\right\vert =s_{i}$ for $1\leq i\leq r.$

\begin{theorem}
\label{th4}Let $r\geq2,$ let $U_{1},\ldots,U_{r}$ be sets of size $n$ and
$M\subset U_{1}\times\cdots\times U_{r}$ satisfy $\left\vert M\right\vert
\geq\alpha n^{r}.$ If
\[
2^{r}\exp\left(  -\frac{1}{r}\left(  \ln n\right)  ^{1/r}\right)  \leq
\alpha\leq1
\]
and the integers $s_{1},\ldots,s_{r}$ satisfy $1\leq s_{1}\cdots s_{r}\leq\ln
n,$ then%
\[
\left\vert B_{M}\left(  s_{1},\ldots,s_{r}\right)  \right\vert \geq\left(
\frac{\alpha}{2^{r}}\right)  ^{rs_{1}\cdots s_{r}}\binom{n}{s_{1}}\cdots
\binom{n}{s_{r}}.
\]

\end{theorem}

\subsubsection*{Remarks}

\begin{itemize}
\item[-] The relations between $\alpha$ and $n$ in the above theorems need
some explanation. First, for fixed $\alpha,$ they show how large must be $n$
to get valid conclusions. But, in fact, the relations are subtler, for
$\alpha$ itself may depend on $n,$ e.g., letting $\alpha=\ln\ln n,$ the
conclusions are meaningful for sufficiently large $n.$

\item[-] Note that, in Theorems \ref{th1}-\ref{th3}, if the conclusion holds
for some $\alpha,$ it holds also for $0<\alpha^{\prime}<\alpha,$ provided $n$
is sufficiently large.

\item[-] As Erd\H{o}s showed in \cite{Erd64}, most graphs with $n$ vertices
and $\left(  1-\varepsilon\right)  \binom{n}{r}$ edges have no $K_{r}\left(
s,\ldots,s\right)  $ for $s\geq c\left(  \log n\right)  ^{1/\left(
r-1\right)  }$ and sufficiently large constant $c=c\left(  \varepsilon\right)
$, independent of $n$. Hence, Theorems \ref{th1}-\ref{th3} are essentially
best possible at least for fixed $\alpha$.

\item[-] Finally, observe that different relations hold for $r=2$, e.g., the
following version of Lemma 2 in \cite{Nik07} corresponds to Theorem \ref{th3}:

\emph{Let }$\left(  \ln n\right)  ^{-1/2}\leq\alpha<1/2,$\emph{ and let }%
$G$\emph{ be a bipartite }$2$\emph{-graph with parts of size }$n$\emph{ with
at least }$\alpha n^{2}$\emph{ edges. Then }$G$\emph{ contains a }%
$K_{2}\left(  s,t\right)  $\emph{ with }$s=\left\lfloor \alpha^{2}\ln
n\right\rfloor $\emph{ and }$t>n^{1-\alpha}$\emph{.}
\end{itemize}

\subsubsection*{Proofs}

First, some definitions.

Suppose $U_{1},\ldots,U_{r}$ are nonempty sets and $M\subset U_{1}\times
\cdots\times U_{r}$; let the integers $s_{1},\ldots,s_{r}$ satisfy
$0<s_{i}\leq\left\vert U_{i}\right\vert ,$ $1\leq i\leq r.$

Define $M^{\prime}\subset U_{1}\times\cdots\times U_{r-1}$ as$\smallskip$

$\qquad M^{\prime}=\left\{  \left(  u_{1},\ldots,u_{r-1}\right)  :\text{there
exists }u\in U_{r}\text{ such that }\left(  u_{1},\ldots,u_{r-1},u\right)  \in
M\right\}  .\smallskip$

For every $R\in B_{M^{\prime}}\left(  s_{1},\ldots,s_{r-1}\right)  ,$
let$\smallskip$

$\qquad N_{M}\left(  R\right)  =\left\{  u:u\in U_{r}\text{ and }\left(
u_{1},\ldots,u_{r-1},u\right)  \in M\text{ for every }\left(  u_{1}%
,\ldots,u_{r-1}\right)  \in R\right\}  ,\smallskip$

$\qquad d_{M}\left(  R\right)  =\left\vert N_{M}\left(  R\right)  \right\vert
.\smallskip$

For every $v\in U_{r},$ let$\smallskip$

$\qquad N_{M}\left(  v\right)  =\left\{  \left(  u_{1},\ldots,u_{r-1}\right)
:\left(  u_{1},\ldots,u_{r-1},v\right)  \in M\right\}  ,\smallskip$

$\qquad d_{M}\left(  v\right)  =\left\vert N_{M}\left(  v\right)  \right\vert
,\smallskip$

$\qquad D_{M}\left(  v\right)  =\left\vert \left\{  R:R\in B_{M^{\prime}%
}\left(  s_{1},\ldots,s_{r-1}\right)  \text{ and }v\in N_{M}\left(  R\right)
\right\}  \right\vert .\smallskip$

Finally, for every integer $s\geq1,$ let\emph{ }%
\[
g_{s}\left(  x\right)  =\left\{
\begin{array}
[c]{cc}%
\binom{x}{s} & \text{if }x>s-1;\\
0 & \text{if }x\leq s-1.
\end{array}
\right.
\]

\begin{proof}
[\textbf{Proof of Theorem \ref{th4}}]We use induction on $r.$ Let first $r=2,$
and by symmetry assume that $s_{1}\geq s_{2}$. Since $g_{s_{2}}\left(
x\right)  $ is convex, we see that
\begin{align*}
\left\vert B_{M}\left(  s_{1},s_{2}\right)  \right\vert  &  =%
{\textstyle\sum\limits_{R\subset U_{1},\left\vert R\right\vert =s_{1}}}
\binom{d_{M}\left(  R\right)  }{s_{2}}=%
{\textstyle\sum\limits_{R\subset U_{1},\left\vert R\right\vert =s_{1}}}
g_{s_{2}}\left(  d_{M}\left(  R\right)  \right) \\
&  \geq\binom{n}{s_{1}}g_{s_{2}}\left(  \binom{n}{s_{1}}^{-1}%
{\textstyle\sum\limits_{R\subset U_{1},\left\vert R\right\vert =s_{1}}}
d_{M}\left(  R\right)  \right)
\end{align*}
On the other hand, we have
\begin{align*}%
{\textstyle\sum\limits_{R\subset U_{1},\left\vert R\right\vert =s_{1}}}
d_{M}\left(  R\right)   &  =%
{\textstyle\sum\limits_{u\in U_{2}}}
\binom{d_{M}\left(  u\right)  }{s_{1}}=%
{\textstyle\sum\limits_{u\in U_{2}}}
g_{s_{1}}\left(  d_{M}\left(  u\right)  \right)  \geq ng_{s_{1}}\left(
\frac{1}{n}%
{\textstyle\sum\limits_{u\in U_{2}}}
d_{M}\left(  u\right)  \right) \\
&  \geq n\binom{\left\vert M\right\vert /n}{s_{1}}\geq n\binom{\alpha n}%
{s_{1}}.
\end{align*}

We have%
\[
\alpha n\geq4\exp\left(  \ln n-\frac{1}{2}\left(  \ln n\right)  ^{1/2}\right)
>2\exp\left(  \frac{1}{2}\ln n\right)  \geq2\ln n.
\]
and so, $\alpha n>2s_{1}.$ Therefore,
\[
n\binom{\alpha n}{s_{1}}\geq n\left(  \frac{\alpha}{2}\right)  ^{s_{1}}%
\binom{n}{s_{1}},
\]
and, since $g_{s_{2}}\left(  x\right)  $ is non-decreasing, we obtain%
\[
\left\vert B_{M}\left(  s_{1},s_{2}\right)  \right\vert \geq\binom{n}{s_{1}%
}g_{s_{2}}\left(  n\binom{n}{s_{1}}^{-1}\binom{\alpha n}{s_{1}}\right)
\geq\binom{n}{s_{1}}g_{s_{2}}\left(  \left(  \frac{\alpha}{2}\right)  ^{s_{1}%
}n\right)  .
\]

Likewise, from%
\[
-\frac{1}{2}\left(  \ln n\right)  ^{1/2}\leq\ln\frac{\alpha}{4}\leq\ln\frac
{1}{4},
\]
we see that $n\geq e^{\left(  \ln16\right)  ^{2}},$ and so,%
\[
\left(  \alpha/2\right)  ^{s_{1}}n\geq\left(  \alpha/2\right)  ^{\ln
n}n=n^{1+\ln\alpha/2}\geq n^{0.3}\geq2\sqrt{\ln n}\geq2s_{2}.
\]
This inequality implies that
\begin{align*}
\left\vert B_{M}\left(  s_{1},s_{2}\right)  \right\vert  &  \geq\binom
{n}{s_{1}}\binom{\left(  \alpha/2\right)  ^{s_{1}}n}{s_{2}}\geq\alpha
^{s_{1}s_{2}}2^{-s_{1}s_{2}-s_{2}}\binom{n}{s_{1}}\binom{n}{s_{2}}\\
&  >\left(  \frac{\alpha}{4}\right)  ^{s_{1}s_{2}}\binom{n}{s_{1}}\binom
{n}{s_{2}},
\end{align*}
completing the proof for $r=2$.

Assume now the assertion true for $r-1;$ we shall prove it for $r.$ We first
show that there exist $W\subset U_{r}$ and $L\subset M$ with $\left\vert
L\right\vert >\left(  \alpha/2\right)  n^{r}$ such that $d_{L}\left(
u\right)  \geq\left(  \alpha/2\right)  n^{r-1}$ for all $u\in W.$ Indeed,
apply the following procedure:

\textbf{Let} $W=U_{r}$,\textbf{ }$L=M;$

\textbf{While }\emph{there exists an }$u\in W$\emph{ with }$d_{L}\left(
u\right)  <\left(  \alpha/2\right)  n^{r-1}$ \textbf{do}

\qquad\emph{Remove }$u$ \emph{from }$W$ \emph{and remove all }$r$\emph{-tuples
containing }$u$ \emph{from }$L.$

When this procedure stops, we have $d_{L}\left(  u\right)  \geq\left(
\alpha/2\right)  n^{r-1}$ for all $u\in W.$ In addition,
\[
\left\vert M\right\vert -\left\vert L\right\vert <\left(  \alpha/2\right)
n^{r-1}n\leq\left(  \alpha/2\right)  n^{r},
\]
implying that $\left\vert L\right\vert \geq\left(  \alpha/2\right)  n^{r},$ as claimed.

Since $g_{s_{r}}\left(  x\right)  $ is convex, and%
\[%
{\textstyle\sum\limits_{R\in B_{L^{\prime}}\left(  s_{1},\ldots,s_{r-1}%
\right)  }}
d_{L}\left(  R\right)  =\left\vert L\right\vert =%
{\textstyle\sum\limits_{u\in W}}
D_{L}\left(  u\right)  ,
\]
we see that
\begin{align}
\left\vert B_{L}\left(  s_{1},\ldots,s_{r}\right)  \right\vert  &  \geq%
{\textstyle\sum\limits_{R\in B_{L^{\prime}}\left(  s_{1},\ldots,s_{r-1}%
\right)  }}
\binom{d_{L}\left(  R\right)  }{s_{r}}=%
{\textstyle\sum\limits_{R\in B_{L^{\prime}}\left(  s_{1},\ldots,s_{r-1}%
\right)  }}
g_{s_{r}}\left(  d_{L}\left(  R\right)  \right) \nonumber\\
&  \geq\left\vert B_{L^{\prime}}\left(  s_{1},\ldots,s_{r-1}\right)
\right\vert g_{s_{r}}\left(  \frac{%
{\textstyle\sum\limits_{R\in B_{L^{\prime}}\left(  s_{1},\ldots,s_{r-1}%
\right)  }}
d_{L}\left(  R\right)  }{\left\vert B_{L^{\prime}}\left(  s_{1},\ldots
,s_{r-1}\right)  \right\vert }\right) \nonumber\\
&  =\left\vert B_{L^{\prime}}\left(  s_{1},\ldots,s_{r-1}\right)  \right\vert
g_{s_{r}}\left(  \frac{%
{\textstyle\sum\limits_{u\in W}}
D_{L}\left(  u\right)  }{\left\vert B_{L^{\prime}}\left(  s_{1},\ldots
,s_{r-1}\right)  \right\vert }\right)  \label{in1}%
\end{align}

On the other hand $s_{1}\cdots s_{r-1}\leq s_{1}\cdots s_{r}\leq\ln n.$ Also,
for every $u\in W,$ we have
\[
\frac{d_{L}\left(  u\right)  }{n^{r-1}}\geq\frac{\alpha}{2};
\]
hence, in view of%
\[
\frac{\alpha}{2}\geq2^{r-1}e^{-\sqrt[r]{\ln n}/r}>2^{r-1}e^{-\sqrt[r-1]{\ln
n}/\left(  r-1\right)  },
\]
we can apply the induction hypothesis to the sets $U_{1},\ldots,U_{r-1},$ the
numbers $s_{1},\ldots,s_{r-1},$ and the set $N_{L}\left(  u\right)  \subset
U_{1}\times\cdots\times U_{r-1}.$ We obtain%
\[
D_{L}\left(  u\right)  \geq\left(  \frac{\alpha/2}{2^{r-1}}\right)  ^{\left(
r-1\right)  s_{1}\cdots s_{r-1}}\binom{n}{s_{1}}\cdots\binom{n}{s_{r-1}}%
\]
for every $u\in W.$ This, together with $\left\vert W\right\vert
\geq\left\vert L\right\vert /n^{r-1}\geq\alpha n/2,$ gives%
\[%
{\textstyle\sum\limits_{u\in W}}
D_{L}\left(  u\right)  \geq\frac{\alpha n}{2}\left(  \frac{\alpha}{2^{r}%
}\right)  ^{\left(  r-1\right)  s_{1}\cdots s_{r-1}}\binom{n}{s_{1}}%
\cdots\binom{n}{s_{r-1}}.
\]

Note that the function $g_{s_{r}}\left(  x/k\right)  k$ is non-increasing in
$k$ for $k\geq1.$ Hence, from
\[
\left\vert B_{L^{\prime}}\left(  s_{1},\ldots,s_{r-1}\right)  \right\vert
\leq\binom{n}{s_{1}}\cdots\binom{n}{s_{r-1}}%
\]
and (\ref{in1}), we obtain
\begin{align}
\left\vert B_{L}\left(  s_{1},\ldots,s_{r}\right)  \right\vert  &  \geq
\binom{n}{s_{1}}\cdots\binom{n}{s_{r-1}}g_{s_{r}}\left(  \binom{n}{s_{1}}%
^{-1}\cdots\binom{n}{s_{r-1}}^{-1}%
{\textstyle\sum\limits_{u\in W}}
D_{L}\left(  u\right)  \right) \nonumber\\
&  \geq\binom{n}{s_{1}}\cdots\binom{n}{s_{r-1}}g_{s_{r}}\left(  \frac{\alpha
}{2}\left(  \frac{\alpha}{2^{r}}\right)  ^{\left(  r-1\right)  s_{1}\cdots
s_{r-1}}n\right)  . \label{in2}%
\end{align}
To continue we need the following

\begin{claim}
The condition
\[
2^{r}\exp\left(  -\frac{1}{r}\left(  \ln n\right)  ^{1/r}\right)  \leq
\alpha\leq1
\]
implies that
\[
\frac{\alpha}{2}\left(  \frac{\alpha}{2^{r}}\right)  ^{\left(  r-1\right)
s_{1}\cdots s_{r-1}}n\geq2s_{r}.
\]

\end{claim}

\emph{Proof} We have%
\begin{equation}
\frac{\alpha}{2}\left(  \frac{\alpha}{2^{r}}\right)  ^{s_{1}\cdots s_{r}}%
\geq\frac{\alpha}{2}\left(  \frac{\alpha}{2^{r}}\right)  ^{\ln n}%
>e^{-\sqrt[r]{\ln n}/r}\left(  e^{-\sqrt[r]{\ln n}/r}\right)  ^{\ln n}=\left(
en\right)  ^{-\sqrt[r]{\ln n}/r}. \label{in3}%
\end{equation}

On the other hand
\[
e^{\ln4-\sqrt[2]{\ln n}/2}=2^{2}e^{-\sqrt[r]{\ln n}/2}\leq2^{r}e^{-\sqrt[r]%
{\ln n}/2},
\]
and so, $\ln n\geq\left(  \ln16\right)  ^{2},$ implying in turn that $n\geq
e^{\left(  \ln16\right)  ^{2}}=16^{\ln16}.$ Routine calculus shows that
$n^{1/2}-4\ln n$ increases for $n\geq16^{\ln16},$ and so,
\[
n^{1/2}-4\ln n\geq\left(  16^{\ln16}\right)  ^{1/2}-4\ln16>0.
\]
Now, from (\ref{in3}) we obtain
\begin{align*}
\frac{\alpha}{2}\left(  \frac{\alpha}{2^{r}}\right)  ^{s_{1}\cdots s_{r}}  &
\geq\left(  en\right)  ^{-\sqrt[r]{\ln n}/r}>\left(  4n\right)  ^{-\left(  \ln
n\right)  /2}=\left(  \frac{1}{2n^{1/2}}\right)  ^{\ln n}\\
&  >\left(  \frac{2\ln n}{n}\right)  ^{\ln n}\geq\left(  \frac{2s_{r}}%
{n}\right)  ^{s_{r}},
\end{align*}
completing the proof of the claim.$\hfill\square$

From (\ref{in2}) and the definition of $g_{s_{r}}\left(  x\right)  $ we see
that
\begin{align*}
\left\vert B_{L}\left(  s_{1},\ldots,s_{r}\right)  \right\vert  &  \geq
\binom{n}{s_{1}}\cdots\binom{n}{s_{r-1}}\left(  \frac{\alpha}{2}\left(
\frac{\alpha}{2^{r}}\right)  ^{\left(  r-1\right)  s_{1}\cdots s_{r-1}%
}\right)  ^{s_{r}}\binom{n}{s_{r}}\\
&  \geq\left(  \frac{\alpha}{2}\right)  ^{s_{r}}\left(  \frac{\alpha}{2^{r}%
}\right)  ^{\left(  r-1\right)  s_{1}\cdots s_{r}}\binom{n}{s_{1}}\cdots
\binom{n}{s_{r}}\\
&  >\left(  \frac{\alpha}{2^{r}}\right)  ^{rs_{1}\cdots s_{r}}\binom{n}{s_{1}%
}\cdots\binom{n}{s_{r}},
\end{align*}
completing the proof.
\end{proof}

\begin{proof}
[\textbf{Proof of Theorem \ref{th3}}]As in the proof of Theorem \ref{th4} we
find $W\subset U_{r}$ and $L\subset M$ with $\left\vert L\right\vert >\left(
\alpha/2\right)  n^{r}$ such that $d_{L}\left(  u\right)  \geq\left(
\alpha/2\right)  n^{r-1}$ for all $u\in W.$ Let
\[
t=\max\left\{  d_{L}\left(  R\right)  :R\in B_{L^{\prime}}\left(  s_{1}%
,\ldots,s_{r-1}\right)  \right\}  .
\]
We have
\[
\left\vert B_{L^{\prime}}\left(  s_{1},\ldots,s_{r-1}\right)  \right\vert
\leq\binom{n}{s_{1}}\cdots\binom{n}{s_{r-1}},
\]
and so
\begin{equation}
t\binom{n}{s_{1}}\cdots\binom{n}{s_{r-1}}\geq t\left\vert B_{L^{\prime}%
}\left(  s_{1},\ldots,s_{r-1}\right)  \right\vert \geq\left\vert L\right\vert
=%
{\textstyle\sum\limits_{u\in W}}
D_{L}\left(  u\right)  . \label{in4}%
\end{equation}
To continue we need the following

\begin{claim}
\label{cl2}The condition $\left(  \ln n\right)  ^{-1/\left(  r-1\right)  }%
\leq\alpha\leq r^{-3}$ implies that
\[
2^{r-1}\exp\left(  -\frac{1}{r-1}\left(  \ln n\right)  ^{1/\left(  r-1\right)
}\right)  \leq\frac{\alpha}{2}\leq1
\]

\end{claim}

\emph{Proof }The upper bound is obvious, so we have to prove that
\[
\ln\frac{\alpha}{2^{r}}\geq-\frac{1}{r-1}\left(  \ln n\right)  ^{1/\left(
r-1\right)  }.
\]
The function $x^{x}$ decreases for $0<x<e^{-1},$ and $\alpha\leq r^{-3};$
hence
\begin{equation}
\alpha\ln\frac{\alpha}{2^{r}}\geq\frac{1}{r^{3}}\ln\frac{1}{r^{3}2^{r}}%
=-\frac{1}{r^{3}}\left(  3\ln r+r\ln2\right)  >-\frac{3r}{r^{3}}\geq-\frac
{1}{r}, \label{in5}%
\end{equation}
and so,%
\[
\ln\frac{\alpha}{2^{r}}>-\frac{1}{\left(  r-1\right)  \alpha}\geq-\frac
{1}{r-1}\left(  \ln n\right)  ^{-1/\left(  r-1\right)  },
\]
completing the proof of the claim.$\hfill\square$

Since for every $u\in W$ we have
\[
\frac{d_{L}\left(  u\right)  }{n^{r-1}}\geq\frac{\alpha}{2},
\]
in view of Claim \ref{cl2}, we may apply Theorem \ref{th4} to the sets
$U_{1},\ldots,U_{r-1},$ the numbers $s_{1},\ldots,s_{r-1},$ and the set
$N_{L}\left(  u\right)  \subset U_{1}\times\cdots\times U_{r-1},$ thus
obtaining%
\[
D_{L}\left(  u\right)  \geq\left(  \frac{\alpha/2}{2^{r-1}}\right)  ^{\left(
r-1\right)  s_{1}\cdots s_{r-1}}\binom{n}{s_{1}}\cdots\binom{n}{s_{r-1}}%
\]
for every $u\in W.$ This, together with $\left\vert W\right\vert
\geq\left\vert L\right\vert /n^{r-1}\geq\alpha n/2,$ gives%
\[%
{\textstyle\sum\limits_{u\in W}}
D_{L}\left(  u\right)  \geq\frac{\alpha n}{2}\left(  \frac{\alpha}{2^{r}%
}\right)  ^{\left(  r-1\right)  s_{1}\cdots s_{r-1}}\binom{n}{s_{1}}%
\cdots\binom{n}{s_{r-1}}.
\]
Substituting this bound in (\ref{in4}), we find that%
\[
t\geq\frac{\alpha}{2}\left(  \frac{\alpha}{2^{r}}\right)  ^{\left(
r-1\right)  s_{1}\cdots s_{r-1}}n\geq\frac{\alpha}{2}\left(  \frac{\alpha
}{2^{r}}\right)  ^{\left(  r-1\right)  \alpha^{r-1}\ln n}n>\left(
\frac{\alpha}{2^{r}}\right)  ^{r\alpha^{r-1}\ln n}n.
\]
Finally, (\ref{in5}) gives
\[
\left(  \frac{\alpha}{2^{r}}\right)  ^{r\alpha^{r-1}\ln n}>e^{-\alpha^{r-2}\ln
n}=n^{-\alpha^{r-2}},
\]
completing the proof of Theorem \ref{th3}.
\end{proof}

\begin{proof}
[\textbf{Proof of Theorem \ref{th2}}]Suppose $r,\alpha,n,$ and $G$ satisfy the
conditions of the theorem. Let $U_{1},\ldots,U_{r}$ be $r$ copies of the
vertex set $V$ of $G,$ and let $M\subset U_{1}\times\cdots\times U_{r}$ be the
set of $r$-vectors $\left(  u_{1},\ldots,u_{r}\right)  $ such that $\left\{
u_{1},\ldots,u_{r}\right\}  $ is an edge of $G.$ Clearly, $\left\vert
M\right\vert \geq r!\left(  \alpha n^{r}/r!\right)  =\alpha n^{r}.$ Theorem
\ref{th3} implies that there exists a set $V_{1}\times\cdots\times
V_{r}\subset M$ such that $V_{i}\subset V$ and $\left\vert V_{i}\right\vert
=s_{i}$ for $1\leq i<r,$ and $\left\vert V_{r}\right\vert >n^{1-\alpha^{r-2}%
}.$ Note that the sets $V_{1},\ldots,V_{r}$ are disjoint, for the edges of $G$
consist of distinct vertices. Hence $V_{1},\ldots,V_{r}$ are the vertex
classes of an $r$-partite subgraph of $G$ with the desired size.
\end{proof}

\end{document}